
\tolerance=10000
\magnification=1200
\raggedbottom

\baselineskip=15pt
\parskip=1\jot

\def\sk{\vskip 3\jot}

\def\heading#1{\vskip3\jot{\noindent\bf #1}}
\def\label#1{{\noindent\it #1}}
\def\QED{\hbox{\rlap{$\sqcap$}$\sqcup$}}


\def\ref#1;#2;#3;#4;#5.{\item{[#1]} #2,#3,{\it #4},#5.}
\def\refinbook#1;#2;#3;#4;#5;#6.{\item{[#1]} #2, #3, #4, {\it #5},#6.} 
\def\refbook#1;#2;#3;#4.{\item{[#1]} #2,{\it #3},#4.}


\def\({\bigl(}
\def\){\bigr)}


\def\al{\alpha}

\def\ep{\varepsilon}
\def\ze{\zeta}
\def\et{\eta}
\def\th{\vartheta}

\def\rh{\varrho}
\def\si{\sigma}

\def\ph{\phi}

\def\ps{\psi}

\def\Ga{\Gamma}
\def\De{\Delta}

\def\La{\Lambda}

\def\Ph{\Phi}

\def\Me{\Omega}



\def\ep{\varepsilon}
\def\ph{\varphi}
\def\th{\vartheta}

\def\mypmod#1{\;({\rm mod}\;#1)}

\def\Ex{{\rm Ex}}
\def\Var{{\rm Var}}

\def\abs#1{\vert#1\vert}
\def\Res{\mathop{\rm Res}}

\def\idle{{\rm\ idle}}

\def\qst{q_*}

\def\sigi{J_{\rm i}}
\def\sigii{J_{\rm ii}}
\def\sigiii{J_{\rm iii}}
\def\sigiv{J_{\rm iv}}

{
\pageno=0
\nopagenumbers
\rightline{\tt spider.web.tex}
\vskip1in

\centerline{\bf The Linking Probability of Deep %
Spider-Web Networks}
\vskip0.5in

\centerline{Nicholas Pippenger}
\centerline{\tt njp@princeton.edu}
\vskip0.25in

\centerline{Department of Computer Science}
\centerline{Princeton University}
\centerline{35 Olden Street}
\centerline{Princeton, NJ 08540}
\vskip0.5in

\noindent{\bf Abstract:}
We consider crossbar switching networks with base $b$ (that is,
constructed from
$b\times b$ crossbar switches), scale $k$ (that is, with $b^k$ inputs,
$b^k$ outputs and $b^k$ links between each consecutive pair of stages)
and depth $l$ (that is, with $l$ stages).
We assume that the crossbars are interconnected according to the 
spider-web pattern, whereby two diverging paths reconverge only after
at least $k$ stages.
We assume that each vertex is independently idle with probability $q$,
the vacancy probability.
We assume that $b\ge 2$ and the vacancy probability $q$ are fixed, and
that
$k$ and $l = ck$ tend to infinity with ratio a
fixed constant $c>1$.
We consider the linking probability $Q$ (the probability that there
exists at least one idle path between a given idle input and a given idle
output).
In a previous paper it was shown that if $c\le 2$, then the
linking probability $Q$ tends to $0$ if $0<q<q_c$ (where $q_c =
1/b^{(c-1)/c}$ is the critical vacancy probability), and tends to
$(1-\xi)^2$ (where $\xi$ is the unique solution of the equation
$\(1-q (1-x)\)^b=x$ in the range $0<x<1$) if $q_c<q<1$.
In this paper we extend this result to all rational $c>1$.
This is done by using generating functions and complex-variable
techniques to estimate the second moments of various random variables
involved in the analysis of the networks.
\vfill\eject
}

\heading{1. Introduction}

We deal in this paper with linking in crossbar switching networks,
a phenomenon not dissimilar to that of percolation in lattices
(as introduced by Broadbent and Hammersley [B] and surveyed by
Grimmett [G]).
An important difference, however, is that while percolation can
be studied in finite subgraphs of a single infinite graph modelling
the lattice, there is no single graph that naturally hosts the 
graph modelling crossbars switching networks we are interested in.
Our first order of business will be to describe these graphs.

A {\it crossbar graph\/}
will be characterized firstly by three parameters:
its {\it base\/} $b\ge 2$, its {\it scale\/} $k\ge 0$ and its
{\it depth \/} $l\ge 0$.
Its vertices are partitioned into $l+1$ {\it ranks}, each
containing $b^k$ vertices, which are labelled with the 
strings of length $k$ over the alphabet $\{0, \ldots, {b-1}\}$.
The vertices in rank $0$ are called {\it inputs}, those in rank 
$l$ are called {\it outputs}, and those in all other ranks are called
{\it links}.
The edges of the graph are partitioned into $l$ {\it stages},
each containing $b^{k+1}$ edges.
For $1\le m\le l$, the edges of stage $m$ are directed out of vertices
in rank $m-1$ and into vertices in rank $m$.
In a {\it spider-web\/} crossbar graph, which is our main
concern in this paper, there is an edge of stage $m$ from vertex $v$ of
rank
$m-1$ to vertex $w$ of rank $m$ if and only if $v$ and $w$ are labelled
by strings that differ at most in position $j$, where 
$j\equiv m \mypmod{k}$.
The edges of each stage are thus partitioned into
$b^{k-1}$ $b\times b$ complete bipartitite graphs (called {\it crossbars}).
The spider-web crossbar graph with base $b$, scale $k$ and depth $l$
will be denoted $G_{b,k,l}$.
We shall see in Section 2 that if $l\ge k$, there are $b^{l-k}$ paths
from a given input to a given output; if $l<k$, there is at most one path
from a given input to a given output.
Our main interest is in spider-web crossbar graphs with $l\ge k$, since in
these graphs any input can be connected by a path to any output; in our
analysis, however, graphs with $l<k$ will occur as subgraphs, so it will
be necessary to allow this case in some of our intermediate results .

We shall assume that each vertex in the graph $G_{b,k,l}$ is independently
assigned the status {\it idle}, with probability $q$
(called the {\it vacancy\/} probability), or {\it busy},
with the complementary probability $p = 1-q$
(called the {\it occupancy\/} probability).
This random assignment of a status to each vertex in a graph will
be called the {\it state\/} of the graph.
Given an input $v$ and output $w$, let $Q_{v,w}$
(called the {\it linking\/} probability) denote the probability that
there exists a path consisting entirely of idle links from $v$ to $w$.
(In general, the linking probability is define as the {\it conditional\/}
probability that there exists an idle path, {\it given\/} that $v$ and $w$
are  themselves idle, but for the probabilisitic model that we are using
this condition is independent.)
We shall see in Section 2 that if $l\ge k$, the probability $Q_{v,w}$ does not depend
on the choice of the input-output pair $(v,w)$, 
so we shall let $Q$ denote the common value of these probabilities.
The complementary probability $P = 1-Q$ (called the {\it blocking\/}
probability), is the probability that all paths between a given
input-output pair $(v,w)$
are broken by a set of busy links.

In practice, the parameter $p$ represents the amount of traffic being
carried by a crossbar network (which one would like to maximize),
and the parameter $P$ represents the fraction of arriving traffic
lost due to congestion within the network (which one would like to 
minimize).
In analysis, however, it is almost always more convenient to work with
the complementary parameters $q$ and $Q$, so we shall work exclusively
with these parameters in what follows.

In practice, a graph $G_{b,k,l}$
would be fixed, and the linking probability $Q$
would be studied as a function of the vacancy probability $q$.
It is found that $Q$ undergoes a rapid transition from a value
near zero to a significantly positive value as $q$ passes through
a neighborhood of $1/b^{(l-k)/(l-1)}$.
This is easily understood in the following way.

Let the random variable $X_{v,w}$ denote the number of idle paths
from $v$ to $w$.
We shall see in Section 2 that if $l\ge k$,
the distribution of $X_{v,w}$ does not depend
on the choice of the input-output pair $(v,w)$, 
so we shall let $X$ denote a random variable with
this common distribution.
Each of the $b^{l-k}$ paths from $v$ to $w$ contains $l-1$ links, which are all idle with probability
$q^{l-1}$.
Thus we have
$$\Ex[X] = b^{l-k} \, q^{l-1}. \eqno(1.1)$$
Thus as $q$ passes through $1/b^{(l-k)/(l-1)}$, the expected number
of idle paths from $v$ to $w$ (called the {\it specific transparency\/})
goes from an exponentially decreasing to an exponentially increasing
function of $k$ and $l$.
This suggests that if $k$ and $l$ tend to infinity in such a way that
their ratio $c = l/k > 1$
remains fixed, while $b$ and $q$ are also held fixed,
then $Q$ will tend to a limit, and that this limit will 
have a discontinuity as $q$ passes through the critical value
$$q_c = 1/b^{(c-1)/c}.$$
(We note that $1<c<\infty$ implies $1/b<q_c<1$.)
Our goal in this paper is to confirm this conjecture, and to determine
the limiting value of $Q$.

Our first step toward this goal,
taken in Section 2, will be to derive the following estimate for the
second moment $\Ex[X^2]$ of $X$.

\label{Theorem 1.1:}
Let $b\ge 2$ and $1/b<q<1$ be fixed.
Then
$$\Ex[X^2] = \Ex[X] \cdot \left(\left({b-1 \over bq-1}\right)^2
b^{l-k} \, q^{l+1} + 1
+ O(l\, b^{l-2k} \, q^{l})  + O(l \, q^{k})\right)$$
as $k,l\to\infty$ with $l\ge k$
and $\(\log (l+1)\)/(k+1)\to 0$.
(The constants in the $O$-terms may depend on $b$ and $q$, but are independent
of $k$ and $l$.)

We observe that this estimate is enough to establish that
the limiting value (if it exists) of $Q$ for $k\to\infty$ and $l = ck$,
cannot be a continuous function of $Q$ as $q$ passes through $q_c$.
Indeed, from Markov's inequality and (1.1) we have
$$Q = \Pr[X\ge 1] \le \Ex[X] = b^{l-k} \, q^{l-1} \to 0 \eqno(1.2)$$
for $q<q_c$.
On the other hand, (1.1) and Theorem 1.1,
together with the inequality
$$\Pr[X\ge 1] \ge {\Ex[X]^2\over \Ex[X^2]}, \eqno(1.3)$$
 imply
$$\eqalignno{
Q =\Pr[X\ge 1]\ge {\Ex[X]^2\over \Ex[X^2]}
&= {(bq-1)^2\over(b-1)^2 q^2 + (bq-1)^2q}
\left(1 + O\left({l\over b^k}\right)
+ O\left({l\, q^k}\right)\right) \cr
& \cr
&\qquad \to {(bq-1)^2\over(b-1)^2 q^2 + (bq-1)^2q} > 0  &(1.4)\cr}$$
for $q=q_c$.
(To verify (1.3), we consider the distribution of $X$ conditioned
on the event $X\ge 1$.
Since $x^2$ is a convex function of $x$, we have
$$\Ex[X^2 \mid X\ge 1] \ge \Ex[X \mid X\ge 1]^2.$$
Multiplying by $\Pr[X\ge 1]^2$ yields
$$\Ex[X^2] \, \Pr[X\ge 1] = 
\Ex[X^2 \mid X\ge 1] \, \Pr[X\ge 1]^2 \ge 
\Ex[X \mid X\ge 1]^2 \, \Pr[X\ge 1]^2 = \Ex[X]^2,$$
which is equivalent to (1.3).)
The inequalities in (1.2) and (1.4) show that the inferior limit
of $Q$ for $q=q_c$ is strictly greater than the limiting value for
$q<q_c$, as claimed.

The argument of the preceding paragraph also sheds some light on the 
condition $\(\log (l+1)\)/(k+1)\to 0$
in Theorem 1.1.
(This condition involves $k+1$ and $l+1$ rather than $k$ and $l$ simply to
avoid dividing by or taking the logarithm of $0$.)
This condition is not the weakest one
sufficient to give an estimate of the form $\Ex[X^2] = O\(\Ex[X]^2\)$,
but it is clear that some upper bound on the growth of $l$ must be imposed,
for with probability $(1-q)^{b^k}$ all the links in a given rank are busy,
disconnecting all input-output pairs.
Thus if $l/(1-q)^{b^k} \to\infty$ we have $Q\to 0$, contradicting 
the implication of (1.3) when $\Ex[X^2] = O\(\Ex[X]^2\)$.

In Section 3, we shall combine Theorem 1.1 with branching-process
arguments from a previous paper, Pippenger [P3], to establish the 
existence and determine the limiting value of $Q$ for $q>q_c$.

\label{Theorem 1.2:}
Let $b\ge 2$ and $0<q<1$ be fixed,
and let $c>1$ be rational.
Then as $k\to\infty$ with $l = ck$, we have
$$Q\to \cases{
0, &if $0<q<q_c$, \cr
& \cr
(1-\xi)^2, &if $q_c < q < 1$, \cr}$$
where
$\xi$ is the unique solution of the equation 
$x = \(1-q (1-x)\)^b$ in the range $0<x<1$.

Some comment is in order concerning the behavior of $(1-\xi)^2$ as a function
of $q$.
The function $f(x) = \(1-q (1-\xi)\)^b$  is a strictly convex function of $x$
for $0 < q\le 1$, since $f''(x) = b(b-1)q^2\(1-q (1-x)\)^{b-2} > 0$
in this range.
Thus the graph of $f(x)$ can intersect the diagonal at most twice in this
range.
There is one intersection at $x=1$, and the conditions $f(0) = (1-q)^b > 0$
and $f'(1) = bq > 1$ imply that there is at least one intersection in
the range $0<x<1$ when $1/b<q<1$.
Thus there is indeed a unique solution of the equation 
$x = \(1-q (1-x)\)^b$ in the range $0<x<1$ when $1/b<q<1$,
and this latter condition is implied by $q_c < q < 1$.
The degree of this equation can be reduced by one (because of the 
solution $x=1$), and it is easy to see that the resulting equation is
irreducible over the field of rational functions of $q$;
thus $\xi$ is an algebraic function of $q$ of degree
$b-1$.
Since $(1-\xi)^2$ is a polynomial in $\xi$, it is also an algebraic 
function of $q$ of degree $b-1$.
Straightforward analysis shows that $Q\to 1$ as $q\to 1$ with
$1-Q = 1-(1-\xi)^2\sim 1-2(1-q)^b$, which may be interpreted as saying that the
main obstacle to linking when $q\to 1$ is complete occupation
of the $b$ links adjacent to the input in the rank $1$, or of the
$b$ links adjacent to the output in rank $l-1$.
As $q\to 1/b$ from above, we have $(1-\xi)^2\sim (bq-1)^2/{b\choose 2}^2$.

Theorem 1.2 was proved, under the additional restriction $c\le 2$,
by Pippenger [P3], so the additional contribution of the current paper
consists of lifting this restriction.
Nevertheless, the techniques used in the current paper go considerably
beyond those employed in the previous paper: the proof of Theorem 1.1
starts with  a detailed combinatorial examination of the intersections
between paths, then uses complex-variable techniques to determine
the asymptotics of the quantities involved.

Spider-web networks were introduced by Ikeno [I]
(though the name has sometimes been used to refer to a broader
class of networks).
They have several optimality properties among networks constructed from
the same type and number of crossbars.
Takagi [T] show that they have the largest linking probability in 
a large class of crossbar networks called ``rhyming'' networks.
Chung and Hwang [C] showed, surprisingly,
that they are {\it not\/} optimal in the 
larger class of ``balanced'' networks.
But Pippenger [P3] showed that they are {\it asymptotically\/} optimal in 
this class for $1<c\le 2$, and the current paper extends this result
to all $c>1$.

The probability distribution on states that  we use was introduced
by Lee [L1] and Le~Gall [L2, L3].
It is by far the easiest to use for analytical purposes, but it suffers
the defect that the set of busy vertices does not form a set of coherent
paths from inputs to outputs.
Models addressing this defect have been introduced by Koverninski\u{\i}
[K] and Pippenger [P1], and the results in Pippenger [P3] have been
extended to these models in Pippenger [P2].
It seems likely that the results of the the present paper can be
similarly extended.

The current paper is self-contained, except for some estimates concerning
branching processes taken from Pippenger [P3].
We have followed the notation of that paper, except that the base,
which was denote $d$ in that paper, is now denoted $b$ (to free the
symbol $d$ for its traditional use in the calculus).
\sk

\heading{2.  The Second Moment}

Our goal in this section is to prove Theorem 1.1.
We begin with a combinatorial result concerning spider-web graphs.

\label{Lemma 2.1:}
The automorphism group of $G_{b,k,l}$ acts
transitively on the paths from inputs to outputs.

\label{Proof:}
Since an automorphism must permute the vertices of each rank among
themselves, an automorphism $\th$ may be regarded as a sequence
$\th = (\th_0, \ldots, \th_l)$ of permutations, one for each rank.
We shall focus on automorphisms in which  each $\th_m$ (for $0\le m\le l$)
is characterized by a string $\th_{m,1} \cdots \th_{m,k}$ of $k$
digits from the alphabet $\{0, \ldots, b-1\}$, and acts on the vertices
of rank $m$ by carrying the vertex labelled $a_1 \cdots a_k$ to
the vertex labelled $a'_1 \cdots a'_k$, where 
$a'_j \equiv a_j + \th_{m,j} \mypmod{d}$ for $1\le j\le k$.
If, for $1\le m\le l$, the string $\th_{m-1}$ differs from the string
$\th_m$ at most in position $j$, where $j\equiv m\mypmod{k}$, then the
sequence $\th = (\th_0, \ldots, \th_l)$ will constitute an
automorphism.

To show that the automorphisms act transitively on the paths, it will 
suffice to show, for some fixed path $u^*$,
that for every path $u$, there is an automorphism that carries $u^*$
to $u$ (since then the inverse of such an automorphism can be used
to carry any other path $u'$ to $u^*$).
A path $u$ may be regarded as a sequence $u = (u_0, \ldots, u_l)$ of 
vertex labels in which, for $1\le m\le l$, the string $u_{m-1}$
differs from the string $u_m$ at most in position $j$, 
where $j\equiv m\mypmod{k}$.
We shall choose for $u^*$ the path $u^* = (0^k, \ldots, 0^k)$.
Then clearly the automorphism $\th = (\th_0, \ldots, \th_l)$ defined
by $\th_m = u_m$ for $0\le m\le l$ carries $u^*$ to $u$.
\QED

\label{Corollary 2.2:}
If $l\ge k$, the graph $G_{b,k,l}$ contains $b^{l-k}$ paths from 
any given input to any given output; if $l<k$, there is at most one
path from any given input to any given output.

\label{Proof:}
If $l\ge k$, every input-output pair is joined by at least one path,
since every position in the 
strings labelling vertices has an opportunity to change at least once.
Thus by Lemma 2.1, every input is joined by the same number of paths.
Since each of the $b^k$ inputs is the origin of $b^l$ paths to outputs, 
there are a total of $b^{l+k}$ paths joining inputs to outputs,
and thus $b^{l-k}$ paths joining each of the $b^{2k}$ input-output pairs.
If $l<k$, there is a path from input $v$ to output $w$ only if 
the labels of $v$ and $w$ agree in the last $k-l$ positions.
Thus $G_{b,k,l}$ breaks into $b^{k-l}$ disjoint components, each containing
$b^l$ vertices in each rank; there is a unique path joining input $v$ to out[ut
$w$ if they belong to the same component, and no path joining them if they belong
to different components.
\QED

\label{Corollary 2.3:}
If $l\ge k$, the automorphism group of $G_{b,k,l}$ acts
transitively on the input-output pairs.

\label{Proof:}
If $k\ge k$, each input-output pair is joined by a path, so
the corollary follows from Lemma 2.1.
\QED

This corollary, together with the fact that the probability distribution
on states of the graph is invariant under automorphisms of the graph,
justifies our earlier assertion that the linking probability
$Q_{v,w}$ and the distribution of the random variable $X_{v,w}$
are independent of the choice of the input-output pair $(v,w)$ when $l\ge k$.
Henceforth we shall focus attention on the input-output pair
$(v^*, w^*) = (0^k, 0^k)$.
If $l\ge k$, this entails no loss of generality.
When $l<k$, we shall only deal with cases in which
the input and output of interest are joined by a path, 
and in these cases there is again no loss of generality.

Fix $b\ge 2$ and $k\ge 1$.
For $l\ge 0$, let $\ph_l(y)$ denote the generating function for the number
of  paths from the input $v^* = 0^k$ to 
the output $w^* = 0^k$ classified
according to the number of links that have labels different from $0^k$;
that is, the coefficient of $y^m$ in $\ph_l(y)$ is the number of
paths from $v^*$ to $w^*$ that have $l-1-m$ links in common
with the path $u^* = (0^k, \ldots, 0^k)$.
Clearly $\ph_l(y)=1$ for $0\le l\le k$, and $\ph_l(y)$ is a polynomial in $y$ of degree $l-1$ 
if $l\ge k+1$.

We are interested in the polynomials $\ph_l(y)$ 
for various values of $l\ge 0$, with $b$ and $k$ fixed.
To determine them, it will be convenient to work with a graph
$G_{b,k}$ that contains as subgraphs all the graphs $G_{b,k,l}$
for various values of $l$.
For any $m\ge l\ge 0$, $G_{b,k,l}$ may be regarded as the subgraph
comprising the vertices in ranks $0$ through $l$ and
the edges in stages $1$ through $l$ of $G_{b,k,m}$.
Thus we may define the infinite graph
$$G_{b,k} = \bigcup_{l\ge k} G_{b,k,l}$$
as the union (inductive limit) of all these graphs.
The graph $G_{b,k}$ has inputs in rank $0$, but all other vertices
will be referred to links.

For $l\ge 0$, the polynomial $\ph_l(y)$ is the generating function
for the number
of  paths from 
the input $v^* = 0^k$ to the link labelled $0^k$ in rank $l$ classified
according to the number of links that have labels different from $0^k$.

Let 
$$\ps(y,z) = \sum_{l\ge 0} \ph_l(y) \, z^l$$
be the generating function for the polynomials $\ph_l(y)$. 
The key to our estimate for the second moment of $X$ is the following
proposition.

\label{Proposition 2.4:}
We have
$$\ps(y,z) = {1 - byz + (b-1)(yz)^{k+1}\over
(1-z)(1-byz) - (b-1)z(1-y)(yz)^k }.$$

\label{Proof:}
In this proof, we shall employ a concise alternative representation
of a path $u = (u_0, \ldots, u_l)$ of length $l\ge 0$
as a string $t = t_1 \cdots t_{k+l}$ of length $k+l$
over the alphabet $B = \{0, \ldots, b-1\}$.
The first $k$ digits $t_1 \cdots t_k$
of $t$ will be the $k$ digits of the label $u_0$.
For $1\le m\le l$, $t_{k+m}$ will be the digit in position $j$
of $u_m$, where $j\equiv m\mypmod{k}$ (the digit of $u_m$ that might be 
different from that of $u_{m-1}$).
Then for $0\le m\le l$, $u_m$ is the string $t_{m+1} \cdots t_{m+k}$.
In particular, the last $k$ digits of $t$ are the $k$ digits of the
label $u_l$ of the link in rank $l$, and the paths from 
the input $v^* = 0^k$ 
to the link labelled 
$0^k$ in rank $l$
are in one-to-one correspondence with the strings
of length $k+l$ over the alphabet $B$ whose first $k$ digits and
last $k$ digits are $0$s.

Given a path $t = 0^k t_{k+1} \cdots t_{l-k} 0^k$, let us overline
each digit $t_{k+m}$ ($1\le m\le l$)
for which $u_{m-1} \not= 0^k$.
The result is a string over the alphabet $B\cup \overline{B}$,
where $\overline{B} = \{\overline{0}, \ldots, \overline{b-1}\}$
is the set of overlined digits.
For $l\ge 0$,
let the language $K_l\subseteq (B\cup \overline{B})^{k+l}$ comprise the
strings obtained in this way for all paths from the input
$v^* = 0^k$ to the link labelled $0^k$ in rank $l$,
and define $K\subseteq (B\cup\overline{B})^{*}$ by
$$K = \bigcup_{l\ge 0} K_l.$$
Then $\ps(y,z)$ is the power series in
$y$ and
$z$ in which the coefficient of  $y^j z^l$ is the number of
strings of length $k+l$ in $K$ in which $j$ digits are
overlined.
Let 
$$L = 0^{-k}\,K = \{t\in(B\cup\overline{B})^{*} : 0^k\,t\in K\}$$
be the language obtained from $K$ by deleting the $k$ initial $0$s from 
each string.
Since none of this initial $0$s are overlined, 
$\ps(y,z)$ is the power series in
$y$ and
$z$ in which the coefficient of  $y^j z^l$ is the number of
strings of length $l$ in $L$ in which $j$ digits are
overlined.

Our next step is to write a regular expression for the language $L$.
Define the alphabets $B' = \{1,\ldots,  b-1\}$ and 
$\overline{B}^{\,\prime} = \{\overline{1}, \ldots, \overline{b-1}\}$.
Then $L$ is described by the regular expression
$$\left(\left(\La + 
\left(\overline{B}^{\,\prime}
\left(\La + \overline{0} + \cdots + \overline{0}^{\,k-1}\right)\right)^*
\overline{B}^{\,\prime}\overline{0}^{\,k-1}\right)0\right)^*, \eqno(2.1)$$
where $\La$ denotes the empty string.
To see this, we observe that a string in $L$ can be uniquely parsed into
into zero or more {\it stretches}, each of which ends with an
unoverlined $0$.
A stretch consists of an unoverlined $0$ optionally preceded by an
{\it excursion}.
An excursion consists a {\it final segment\/} preceded by zero or more
{\it preliminary segments}.
A final segment consists of a digit from 
$\overline{B}^{\,\prime}$ followed by exactly
$k-1$ overlined $0$s.
A preliminary segment consists of a digit from 
$\overline{B}^{\,\prime}$ followed by at
most $k-1$ overlined $0$s.
Clearly a final segment is described by the regular expression
$\overline{B}^{\,\prime}\overline{0}^{\,k-1}$, and 
a preliminary segment is described by the regular expression
$\overline{B}^{\,\prime}
\(\La + \overline{0} + \cdots + \overline{0}^{\,k-1}\)$
Thus an excursion is described by the regular expression
$$\left(\overline{B}^{\,\prime}
\left(\La + \overline{0} + \cdots + \overline{0}^{\,k-1}\right)\right)^*
\overline{B}^{\,\prime}\overline{0}^{\,k-1},$$
and a stretch is described by the regular expression
$$\left(\La + 
\left(\overline{B}^{\,\prime}
\left(\La + \overline{0} + \cdots + \overline{0}^{\,k-1}\right)\right)^*
\overline{B}^{\,\prime}\overline{0}^{\,k-1}\right)0.$$
Thus the strings in $L$ are described by the regular expression (2.1).

We now observe that the regular expression (2.1) is {\it unambiguous\/}
in the following sense: a string described by a subexpression
$R+S$ is described by $R$ or by $S$ (but not both), a string
$t$ described by a subexpression $RS$ has unique parsing $t=rs$
such that $r$ is described by $R$ and $s$ is described by $S$,
and a string $t$ described by a subexpression $S^*$ has a unique
parsing $s = s_1 \cdots s_n$ with $n\ge 0$ such that
$s_1, \ldots, s_n$ are each described by $S$.

For an unambiguous regular expression, if $\ps_R(y,z)$ and $\ps_S(y,z)$
are the generating functions counting the strings described by
subexpressions $R$ and $S$, respectively, then 
$\ps_R(y,z) + \ps_S(y,z)$, $\ps_R(y,z) \, \ps_S(y,z)$ and
$1/\(1 - \ps_S(y,z)\)$ 
are the generating functions counting the strings described by
the subexpressions $R + S$, $RS$ and $S^*$, respectively.

Thus the final segments are counted by the generating function
$(b-1)(yz)^k$ and the preliminary  segments 
are counted by the generating function 
$$(b-1)yz\(1 + yz + \cdots + (yz)^{k-1}\) 
= {(b-1)(yz - (yz)^{k+1})\over1 - yz}.$$
The excursions are counted by 
$${(b-1)(yz)^k\over \displaystyle 1 -  {(b-1)(yz - (yz)^{k+1})\over
1 - yz}} 
= { (b-1)((yz)^k - (yz)^{k+1}) 
\over 1 - yz - (b-1)(yz - (yz)^{k+1})},$$
and the stretches are counted by
$$\left( 1 + { (b-1)((yz)^k - (yz)^{k+1}) 
\over 1 - yz - (b-1)(yz - (yz)^{k+1})}\right)z
= {  z - yz^2- (b-1)z(yz - (yz)^k)
\over 1 - yz - (b-1)(yz - (yz)^{k+1})}.$$
Thus the strings in $L$ are counted by
$${1\over 1 - \displaystyle
{  z - yz^2- (b-1)z(yz - (yz)^k)
\over 1 - yz - (b-1)(yz - (yz)^{k+1})}}
= {1 - byz + (b-1)(yz)^{k+1}\over
(1-z)(1-byz) - (b-1)z(1-y)(yz)^k },$$
which completes the proof of the proposition.
\QED

\label{Proposition 2.5:}
Let $b\ge 2$ and $0<q<1$ be fixed.
Then as $k\to \infty$ 
and $l\ge 0$ behaves in such a way that in such a way that 
$\(\log (l+1)\)/(k+1)\to 0$,
we have
$$\ph_l(q) = \left({b-1 \over bq-1}\right)^2 b^{l-k} q^{l+1} + 1
+ O(l b^{l-2k} q^l) + O(l q^k) + O(l q^l).$$
(The constants in the $O$-terms may depend on $b$ and $q$, but are independent
of $k$ and $l$.)

\label{Proof:}
Write $A(z) = 1 - bqz + (b-1)(qz)^{k+1}$ and 
$B(z) = (1-z)(1-bqz) - (b-1)z(1-q)(qz)^k$, so that
$\ps(q,z) = A(z)/B(z)$.
Then from Cauchy's formula we have
$$\eqalignno{
\ph_l(q)
&= {1\over 2\pi i} \oint_{\Ga_0} {\ps(q,z)\, dz\over z^{l+1}} \cr
&= {1\over 2\pi i} \oint_{\Ga_0} 
{ A(z)\over B(z) }\;\;
{dz\over z^{l+1}}, &(2.2)\cr
}$$
where $\Ga_0$ is a contour taken counterclockwise around a circle
$\abs{z} = \ep$
centered at $0$ and having radius $\ep$ suficiently small
to exclude all other singularities
of the integrand.

To make further progress, we must estimate the locations of these
other singularities, which are poles at the values of $z$ for which
the denominator $B(z)$ vanishes.
One such singularity is at $z = 1/q$.
Let
$$\ze_1 = {1\over q}\left(1 - {1\over l}\right),$$
and let $\Ga_1$ be a contour taken counterclockwise around the circle
$\abs{z} = \ze_1$ centered at $0$ and having radius $\ze_1$.
As $z$ traverses this contour, the magnitude of
the first term $(1-z)(1-bqz)$
of $B(z)$ satisfies the lower bound
$$\eqalign{
\abs{(1-z)(1-bqz)}
&= \abs{1-z} \cdot \abs{1-bqz} \cr
&\ge \left({1\over q} - 1 - {1\over ql}\right)
\left(b - 1 - {b\over l}\right) \cr
&\ge \left({1\over q} - 1\right)
\left(b - 1\right) 
- {bq-1\over ql}, \cr
}$$
since the minimum occurs when $z$ is real and positive.
The magnitude of the second term $(b-1)z(1-q)(qz)^k$, on the other hand,
satisfies the upper bound
$$\eqalign{
\abs{(b-1)z(1-q)(qz)^k}
&=
(b-1)\left({1\over q} - 1\right)\left(1 - {1\over l}\right)^{k+1} \cr
&\le (b-1)\left({1\over q} - 1\right) e^{-k/l} \cr
&\le (b-1)\left({1\over q} - 1\right) 
\left(1 - {(e-1)k\over el}\right). \cr
}$$
Here we have used the inequality $1-x \le e^{-x}$, which holds for all
$x$ because the graph of the convex function $e^{-x}$ lies above that 
of $1-x$, its tangent at $x=0$, and the inequality
$e^{-x} \le 1 - (e-1)x/e$, which holds for $0\le x\le 1$ because
the graph of the convex function $e^{-x}$ lies below that of 
$1 - (e-1)x/e$, its chord across the interval $0\le x\le 1$.
Thus for all sufficiently large $k$ 
(specifically, for $k > (bq-1)e/(b-1)(1-q)(e-1)$), we have the bound
$$\abs{B(z)} = \Me\left({1\over l}\right)$$
for $z$ on the contour $\Ga_1$.
Since we also have $A(z) = O(1)$ for $z$ on $\Ga_1$, we have the estimate
$$\eqalignno{
{1\over 2\pi i} \oint_{\Ga_1} 
{ A(z)\over B(z) }\;\;{dz\over z^{l+1}}
&= O(l\,q^l). &(2.3)\cr}$$

Furthermore, as $z$ traverses the contour $\Ga_1$, the value of the 
first term $(1-z)(1-bqz)$  in $B(z)$ encircles the origin twice, 
since it is a quadratic polynomial.
Since the second term 
$(b-1)z(1-q)(qz)^k$ has strictly smaller magnitude, the value of
$B(z)$ also encircles the origin twice.
It follows that the denominator of $B(z)$ has exactly two zeroes
inside the contour $\Ga_1$.
These are perturbations of the zeros of the first term:
the zero of the first term at $z=1$ is perturbed to one at
$$z = \ze_2 = 1 + O(q^k), \eqno(2.4)$$
and the zero of the first term at $z = 1/bq$ is perturbed to one at
$$z = \ze_3 = {1\over bq}
\left(1 - {(b-1)(1-q)\over (bq-1)\,b^k}
+ O\left({k\over b^{2k}}\right)\right). \eqno(2.5)$$
The condition $\(\log (l+1)\)/(k+1)\to 0$ ensures that the $O$-terms in (2.4) and (2.5)
have smaller orders of magnitude than the terms preceding them.
We observe that $0 < \ze_3 < \ze_2 < \ze_1$,
so that $0$, $\ze_3$ and $\ze_2$ lie inside $\Ga_1$,
and lie in that order along the real axis.
Let $\Ga_2$ be a contour taken counterclockwise around a circle
$\abs{z - \ze_2} = \ep$ centered at 
$\ze_2$ and having radius $\ep$ suficiently small
to exclude all other singularities
of the integrand, and let
$\Ga_3$ be a contour taken counterclockwise around a circle
$\abs{z - \ze_3} = \ep$ centered at 
$\ze_3$ and having radius $\ep$ suficiently small
to exclude all other singularities
of the integrand.
Since the integral of an analytic function around a contour
depends only on the homology class of the contour in the domain
of analyticity of the function, and since $\Ga_0$ is homologous
to $\Ga_1 - \Ga_2 - \Ga_3$
(indeed, $\Ga_0$ is homotopic to a contour that joins a forward traversal of
$\Ga_1$ with reverse traversals of $\Ga_2$ and $\Ga_3$ by
cancelling traversals of segments 
$[\ze_3+\ep, \ze_2 - \ep]$ and $[\ze_2 + \ep, \ze_1]$
of the real axis), from (2.2) we have
$$\eqalignno{
\ph_l(q)
&= {1\over 2\pi i} \oint_{\Ga_1} 
{ A(z)\over B(z) }\;\;{dz\over z^{l+1}} \cr
&\qquad - {1\over 2\pi i} \oint_{\Ga_2} 
{ A(z)\over B(z) }\;\;{dz\over z^{l+1}} \cr
&\qquad - {1\over 2\pi i} \oint_{\Ga_3} 
{ A(z)\over B(z) }\;\;{dz\over z^{l+1}}. &(2.6)\cr
}$$

The first integral in (2.6) has already been estimated in (2.3).
The remaining integrals each encircle just one singularity of the integrand,
and thus they can be evaluated by Cauchy's formula.
If $\ze$ is a simple pole of the integrand, an $\Ga$ is a contour taken
clockwise around just this singularity of the integrand, then we have
$$
\eqalign{{1\over 2\pi i} \oint_{\Ga} 
{ A(z)\over B(z) }\;\;{dz\over z^{l+1}}
&= \Res_{z=\ze} \;\;{ A(z)\over B(z) }\;\;{1\over z^{l+1}} \cr
&= {A(\ze)\over B'(\ze)}\;\;{1\over \ze^{l+1}}. \cr
}$$
For the integral around $\Ga_2$ we have
$A(\ze_2) = -(bq-1) + O(q^k)$ and $B'(\ze_2) = (bq - 1) + O(k\,q^k)$, so that
$$- {1\over 2\pi i} \oint_{\Ga_2} 
{ A(z)\over B(z) }\;\;{dz\over z^{l+1}}
= 1 + O(l\,q^k). \eqno(2.7)$$
For the integral around $\Ga_3$ we have 
$A(\ze_3) = (b-1)^2\/(bq-1)\,b^{k+1} + O(k/b^{2k})$ and
$B'(\ze_3) = -(bq-1) + O(k/b^k)$, so that
$$- {1\over 2\pi i} \oint_{\Ga_3} 
{ A(z)\over B(z) }\;\;{dz\over z^{l+1}}
= \left({b-1 \over bq-1}\right)^2 b^{l-k}\,q^{l+1} 
+ O(l\,b^{l-2k}\,q^l). \eqno(2.8)$$
Substituting the estimates (2.3), (2.7) and (2.8) into (2.6)
completes the proof of the proposition.
\QED

We observe that by extending the asymptotic expansions in (2.4) and (2.6),
it would be possible to extend the expansions in (2.7) and (2.8), and thus
reduce their contibutions to the error terms in Proposition 2.4.
The error term in (2.3), however, cannot be improved without taking account
of the zeroes of $B(z)$ outside the circle $\abs{z} = 1/q$, which will in
general contribute oscillatory terms to the expansion of $\ph_l(q)$.

\label{Proof of Theorem 1.1:}
By Corollary 2.3, we may take $X$ to be the number of idle paths from
$v^* = 0^k$ to $w^* = 0^k$.
We then have
$$\eqalignno{
&= \sum_{u': v^*\to w^*} \; \sum_{u: v^*\to w^*}
\Pr[u \idle, u' \idle] \cr
& = \sum_{u': v^*\to w^*} \;\Pr[u' \idle] \; \sum_{u: v^*\to w^*}
\Pr[u \idle \mid u' \idle], &(2.9)\cr
}$$
where the sums are over all paths from $v^*$ to $w^*$.
For each path $u'$, we can find by Proposition 2.1 an automorphism $\th$ that
carries $u'$ to the the path $u^*$ in which all links are labelled $0^*$.
Applying this automorphism to both $u$ and $u'$ gives
$\Pr[u \idle \mid u' \idle] = 
\Pr[\th(u) \idle \mid u^* \idle]$, since the probability
distribution on states of the graph is invariant under automorphisms.
Furthermore,
$$\eqalign{
 \sum_{u: v^*\to w^*}
\Pr[u \idle \mid u' \idle]
&=
\sum_{u: v^*\to w^*} \Pr[\th(u) \idle \mid u^* \idle] \cr
&= \sum_{u: v^*\to w^*} \Pr[u \idle \mid u^* \idle], \cr
}$$
since both right-hand sides sum the same terms in different orders.
Thus the inner sum in (2.9) is independent of $u'$, and we have
$$\Ex[X^2] = \sum_{u': v^*\to w^*} \;\Pr[u' \idle] \; 
\sum_{u: v^*\to w^*} \Pr[u \idle \mid u^* \idle],$$
so that $\Ex[X^2]$ factors as the product of two sums.
The first sum is just $\Ex[X]$.
To evaluate the second sum, we observe that 
$\Pr[u \idle \mid u^* \idle]$ is just $q^j$, where $j$ is the 
number of links on $u$ that are not labelled $0^k$.
Thus the second sum is $\ph_l(q)$, and the theorem follows from Proposition
2.5.
\QED
\sk

\heading{3. The Linking Probability}

Our goal in this section is to prove Theorem 1.2.
Thus in this section we shall always assume that $b\ge 2$ and $0<q<1$
are fixed, and that $k\to\infty$ and $l=ck$ for some fixed rational
$c>1$.
Thus the constants in $O$-tems may depend on $c$ as well as on $b$ and $q$,
but not otherwise on $k$ or $l$.
We shall also assume that $k$ is even; the case of odd $k$ requires
only that $k/2$ be replaced by $\lfloor k/2\rfloor$ and $\lceil k/2\rceil$
in appropriate ways.

\label{Lemma 3.1:}
Let $G^*_{b,k,l}$ be the graph obtained from $G_{b,k,l}$ by reversing
the direction of its edges and exchanging the roles of its inputs and outputs.
Then $G^*_{b,k,l}$ is isomorphic to $G_{b,k,l}$.

\label{Proof:}
The isomorphism takes the vertex with label $a_1 \cdots a_k$ in rank $m$ of 
$G_{b,k,l}$ to the vertex with label
$a^*_1 \cdots a^*_k$ in rank $l-m$ of 
$G^*_{b,k,l}$, where $a^*_i = a_j$ with $j\equiv l+1-i\mypmod{k}$
(or {\it vice versa}: it is an involution).
\QED

Lemma 3.1 establishes a symmetry between $G_{b,k,l}$ and
$G^*_{b,k,l}$, which we shall invoke by use of the term ``dually''.
(When $l$ is even, $G_{b,k,l}$ is in fact isomorphic to a graph
with manifest bilateral symmetry, as is shown in the Appendix of Pippenger [P3].)

\label{Lemma 3.2:}
Let $\langle G_{b,k,l}\rangle_{m,n}$ , with $0\le m\le n\le l$,
be the subgraph of $G_{b,k,l}$
comprising the vertices in ranks $m$ (now considered inputs)
through $n$ (now considered outputs) and the edges
in stages $m+1$ through $n$.
Then $\langle G_{b,k,l}\rangle_{m,n}$ is isomorphic to $G_{b,k,n-m}$.

\label{Proof:}
The isomorphism takes the vertex with label $a_1 \cdots a_k$ in rank
$h$ of $G_{b,k,n-m}$ to the vertex with label
$a'_1 \cdots a'_k$ in rank $m+h$ of $\langle G_{b,k,l}\rangle_{m,n}$,
where $a'_i = a_j$ with $j\equiv i+m\mypmod{k}$.
\QED

\label{Corollary 3.3:}
Between any given input and any given output of 
$\langle G_{b,k,l}\rangle_{m,n}$, there are $b^{n-m-k}$ paths if
$n-m\ge k$, and either one path or none if $n-m<k$.

\label{Proof:}
Immediate from Lemmas 3.2 and 2.2.
\QED
 
We begin with the upper bound to $Q$.
For $0<q<q_c$, where $q_c = 1/b^{(c-1)/c}$, we have $Q\to 0$ by (1.2).
For $q_c < q < 1$, we shall use the following lemma from Pippenger [P3]
(Corollary 4.2).

\label{Lemma 3.4:}
Let $T_r$ be a complete balanced $b$-ary tree of depth $r$, and let 
each vertex of $T_r$ except for the root be considered idle with probability
$q$ independently.
Let the random variable $Z_r$ denote the number of leaves (vertices at depth
$r$) for which every vertex on the path from the root (exclusive) to the 
leaf (inclusive) is idle.
Then we have
$$\Pr[Z_r = 0] = \xi + O(\et^r)$$
as $r\to\infty$ with $b\ge 2$ and $1/b < q < 1$ fixed, where 
$\xi$ is the unique solution of the equation
$\(1-q (1-\xi)\)^b=\xi$ in the range $0<\xi<1$, and
$\et = b\(1 - q(1-\xi)\)^{b-1} < 1$.

Now set $r = k/2$ and $s=l-k/2$.
The paths from an input $v$ to links in rank $r$ of $G_{b,k,l}$ form a tree
isomorphic to $T_r$ (if we ignore the directions of the edges), and
the paths from links in rank $s$ to an output $w$ for a disjoint tree
isomorphic to $T_r$.
Thus the number of links 
$u$ in rank $r$ for which all the links on the path from
$v$ to $u$ are idle is a random variable $U$
with the same distibution as $Z_r$.
Dually, the number of  links $u$ in rank $s$ for which all the links on the
path from $u$ to $w$ are idle is an independent random variable 
$U'$ with the same
distribution as $Z_r$.
If $v$ and $w$ are linked, then we must have $U\ge 1$ and $U'\ge 1$,
so by Lemma 3.1 we have
$$Q\le \Pr[U\ge 1, U'\ge 1] = (1-\xi)^2 + O(\et^r).$$
This completes the upper bound for Theorem 1.2.
\sk

We now turn to the lower bound for Theorem 1.2.
Since this result has been proved for $c\le 2$ in Pippenger [P3],
we shall assume that $c>2$.
(This assumption could of course be avoided, but it would require
a more complicated choice of parameters and consideration of cases.)
For $0<q<q_c$, there is nothing to prove,
since $Q$ is certainly non-negative.
For $q_c < q < 1$, we shall use the following lemma from Pippenger [P3]
(Lemma 8.1).

\label{Lemma 3.5:}
With  $Z_r$ as in Lemma 3.4 and $1\le H\le (bq)^r$, we have
$$\Pr[Z_r \le H] \le \xi + O\left(\({H / (bq)^r}\)^\al \right)$$
as $r\to\infty$ with $b\ge 2$ and $1/b < q < 1$ fixed, where
$\al = \log(1/\et)\big/\log(bq)$ and $\et$ is as in Lemma 3.4.

Supposing that $q_c<q<1$, we shall define
$$\qst = q_{c-1} \, q^{1/(c-1)^2}.$$
We observe that $q<1$ implies $\qst<q_{c-1}$, and that
$q_c<q$ implies $\qst<q$.

\label{Lemma 3.6:}
Let $k\to\infty $ and $l=ck$, with $b\ge 2$, $q_c<q<1$ and $c>2$ fixed.
Then for all sufficiently large $k$, we have
$$\ps_h(\qst)\le k$$ 
for all $0\le h\le l-k$.

\label{Proof:}
From Proposition 2.4 we have
$$\ph_h(\qst) = \left({b-1 \over b\qst-1}\right)^2 b^{h-k} \qst^{h+1} + 1
+ O(h b^{h-2k} \qst^h) + O(h \qst^k) + O(h \qst^h).$$
Since $\qst<q_{c-1}$ and $h\le l-k$, each term  is O(1), and thus at most
$k$ for all sufficiently large $k$.
\QED

Let
$$H = \lceil (b\qst)^r\rceil.$$
We observe that $v$ and $w$ will be linked if the following three events occur.
\medskip

\item{I.} The input
$v$ is joined by paths containing only idle links to all the links in
a set $V$ containing
at least $H$ idle links in rank $r$.

\item{II.}
All the links in a  set $W$ containing at least $H$ idle links in rank $s$ are
joined y path containing only idle links to the output $w$.

\item{III.}
There is at least one path containing only idle links from some link in $V$
to some link in $W$.
\medskip

By Lemma 3.5, we have
$$\Pr[I] \ge 1 - \xi + O\left(\({\qst / q}\)^r\right),$$
and since $\qst<q$ we have $\Pr[I]\to 1-\xi$.
Dually, we have by Lemma 3.5
$$\Pr[II] \ge 1 - \xi + O\left(\({\qst / q}\)^r\right),$$
and thus also $\Pr[II]\to 1-\xi$.
Since events I and II are independent, we have $\Pr[I,II]\to (1-\xi)^2$.
Thus to complete the proof of the lower bound for Theorem  1.2, it will suffice
to show that
$$\Pr[III\mid I,II] \to 1.$$
Event III depends on events I and II through the sets $V$ and $W$.
We can avoid having to consider this dependence by showing that
$\Pr[III]\to 1$ for {\it any\/} sets $V$ and $W$ each containing at least
$H$ links.
Thus it will suffice to prove the following propostion.

\label{Proposition 3.7:}
Let  $V$ and $W$ be any sets of links in ranks $r$ and $s$, respectively,
each containing at least $H$ links.
Then
$$\Pr[III]\to 1$$
as $k\to\infty$ with  $l = ck$ and $b\ge 2$, $c>2$ and $q_c<q<1$ fixed.

\label{Proof:}
Since $\Pr[III]$ can only increase if links are added to $V$ or $W$,
we may assume that $V$ and $W$ each contain {\it exactly\/} $H$ links.
And since $\Pr[III]$ can only increase if $q$ is increased, it will
suffice to estimate $\Pr[III]$ assuming the vacancy probability to be
$\qst<q$ rather than $q$.

Let the random variable $Y$ be the number of paths containing only idle links
joining some link in $V$ (exclusive) to some link in $W$ (exclusive).
Then event III is equivalent to $Y\ge 1$, so it will suffice to show that
$\Pr[Y=0]\to 0$.
To do this, we shall use Chebyshev's inequality:
$$\Pr[Y=0] \le {\Var[Y]\over \Ex[Y]^2}.$$
Each path from a link in rank $r$ (exclusive) to a link in rank $s$ exclusive
contains $s-r-1 = l - k - 1$ links.
Since each of  these links is independently idle with probability $\qst$,
the probability that such a path contains only idle links is
$\qst^{l-k-1}$.
By Corollary 3.3, the number of such paths joining a given link in rank $r$ with a given
link in rank $s$ is $b^{s-r-k} = b^{l-2k}$.
Since there are $H$ links in each  of $V$ and $W$, we have
$$\eqalign{
\Ex[Y] 
&= H^2\,b^{l-2k} \, \qst^{l-k-1}. \cr
}$$
Next we must estimate $\Var[Y]$.
We have
$$\eqalign{
\Var[Y] 
&= \sum_{u':V\to W} \sum_{u:V\to W}
\(\Pr[u,u'\idle] - \Pr[u\idle]\, \Pr[u'\idle]\) \cr
&= \sum_{u':V\to W} \Pr[u'\idle] \sum_{u:V\to W}
\(\Pr[u\idle \mid u'\idle] - \Pr[u\idle]\) \cr
}$$
Here each sum is over all $H^2$ paths joining a link in $V$ to a link
in $W$, so there are $H^4$ terms in all.
If $u$ is a path from a link in rank $r$ to a link in rank $s$, let
$\rh(u)$ denote the link in rank $r$ and $\si(u)$ the link in rank $s$.
By Proposition 2.1, we may assume (as in the proof of Theorem 1.1)
that $u' = u^*$ is part of a path 
from $v^* = 0^k$ through $\rh(u') = 0^k$ and $\si(u') = 0^*$ to $w^* = 0^k$,
in which all the links  have label $0^k$.
Thus we have
$$\Var[Y] = H^2\,b^{l-2k} \, \qst^{l-k-1} \sum_{u:V\to W}
\(\Pr[u\idle \mid u^*\idle] - \Pr[u\idle]\).$$
The factor $H^2\,b^{l-2k} \, \qst^{l-k-1}$ multiplying the sum is $\Ex[Y]$,
so to show that $\Var[Y]/\Ex[Y]^2 \to 0$, it will suffice to show that
$J/\Ex[Y]\to 0$, where
$$J = \sum_{u:V\to W}
\(\Pr[u\idle \mid u^*\idle] - \Pr[u\idle]\).$$
We now partition the paths $u$ into four classes:
\medskip

\item{i.} those for which $\rh(u) = \si(u) = 0^k$;

\item{ii.} those for which $\rh(u) \not=  0^k$ but $\si(u) = 0^k$;

\item{iii.} those for which $\rh(u) =  0^k$ but $\si(u) \not= 0^k$; and

\item{iv.} those for which $\rh(u) \not=  0^k$ and $\si(u) \not= 0^k$.
\medskip

\noindent We shall denote the contributions to $J$ over these four
classes by $\sigi$, $\sigii$, $\sigiii$ and $\sigiv$, respectively, and 
estimate them in turn.

For $\sigi$, we have
$$\eqalign{
\sigi
&\le \sum_{u:0^k\to 0^k} \Pr[u\idle\mid u^*\idle] \cr
&= \ps_{s-r}(\qst) \cr
&\le k \cr }$$
by Lemma 3.6.
Thus we have
$$\eqalign{
{\sigi\over\Ex[Y]} 
&\le {k\over H^2\,b^{l-2k} \,  \qst^{l-k-1}} \cr
&\le {k\over b^{l-k} \, \qst^{l}} \cr
&\to 0, \cr
}$$
since $\qst>q_c$.

For $\sigii$, we have
$$\sigii \le \sum_{V\setminus\{0^k\}\to 0^k} \Pr[u\idle\mid u^*\idle].$$
To estimate $\Pr[u\idle\mid u^*\idle]$, let $i$ be the first rank for which
a link in $u$ has label $0^k$.
Since there are two distinct paths
in $\langle G_{b,k,l}\rangle_{0,i}$ from $v^*$ through $\rh(u^*)=0^k$
and $\rh(u)\not=0^k$ to this link, we must have $i\ge k+1$ by Corollary 3.3.
Thus we have
$$\eqalign{
\sigii 
&\le (H-1)\left(
\sum_{k+1\le i\le k+r} \qst^{i-r-1} \ps_{s-i}(\qst) + 
\sum_{k+r+1\le i\le s} b^{i-r-k} \qst^{i-r-1} \ps_{s-i}(\qst)
\right) \cr
&\le (H-1)k\left(
\sum_{k+1\le i\le k+r} \qst^{i-r-1}  + 
\sum_{k+r+1\le i\le s} b^{i-r-k} \qst^{i-r-1} 
\right) \cr
}$$
where the factor of $H-1$ accounts for the choice of $\rh(u)\in V\setminus\{0^k\}$,
the factors preceding $\ps_{s-i}(\qst)$ in the
sums account for the probability that all the links
on $u$ between ranks $r$ (exclusive) and $i$ (exclusive) are idle,
the factors of $\ps_{s-i}(\qst)$ account for the probability that all
the links of $u$ between  ranks $i$ and $s$ that are not labelled $0^k$
are idle, and we have bounded $\ps_{s-i}(\qst)$ using Lemma 3.6.
Bounding the sums by the number of terms (at most $s-r = l-k$)
times the largest term (the first for the first sum, and the last for the
second), we have
$$\sigii \le (H-1)k\,(l-k)(\qst^{k/2} + b^{l-2k} \qst^{l-k-1}).$$
Thus we have
$$\eqalign{
{\sigii\over\Ex[Y]}
&\le {k(l-k)(\qst^{k/2} + b^{l-2k} \qst^{l-k-1})
\over H\,b^{l-2k} \,  \qst^{l-k-1}} \cr
&\le k(l-k)
\left({1\over b^{l-3k/2}\,\qst^{l-k}} + {1\over (b\qst)^{k/2}}\right) \cr
&\to 0, \cr
}$$
since $b^{c-3/2}\qst^{c-1} > b^{c-3/2}q_c^{c-1} = b^{-1/2} q_c^{-1}
> b^{-1/2} q_2^{-1} = 1$ (because
$\qst>q_c$, $b^{c-1}q_c^c = 1$, $q_c < q_2$ and $bq_2^2 = 1$)
and $b\qst > 1$ (because $\qst > q_\infty = 1/b$).

Dually, we have $\sigiii/\Ex[Y]\to 0$.

Finally, for $\sigiv$ we have
$$\eqalign{
\sigiv 
&= \sum_{\scriptstyle u:V\setminus\{0^k\}\to W\setminus\{0^k\}}
\(\Pr[u\idle \mid u^*\idle] - \Pr[u\idle]\) \cr
&= \sum_{\scriptstyle u:V\setminus\{0^k\}\to W\setminus\{0^k\}
\atop \scriptstyle u\cap u^*\not=\emptyset}
\(\Pr[u\idle \mid u^*\idle] - \Pr[u\idle]\) \cr
&\le \sum_{\scriptstyle u:V\setminus\{0^k\}\to W\setminus\{0^k\}
\atop \scriptstyle u\cap u^*\not=\emptyset}
\Pr[u\idle \mid u^*\idle], \cr
}$$
since if $u\cap u^*=\emptyset$, the events ``$u\idle$'' and ``$u^*\idle$'' are
independent and the summand $\Pr[u\idle \mid u^*\idle] - \Pr[u\idle]$ vanishes.
Given a path $u$ with $u\cap u^*\not=\emptyset$, let $i$ be  the first
rank in which $u$ has a link with label $0^k$, and let  $j\ge i$ be the last
such rank.
As in case ii, we have $k+1\le i$, and dually we have $j \le l-k-1$.
Thus we have
$$\eqalign{
\sigiv 
&\le 
(H-1)^2 \bigg(\sum_{k+1\le i\le k+r} \;
\sum_{\scriptstyle l-k-r\le j\le l-k-1\atop\scriptstyle i\le j}
\qst^{i-r-1}\ps_{j-i}(\qst)\qst^{s-j-1}  \cr
&\qquad\qquad+ 
\sum_{k+1\le i\le k+r} \; 
\sum_{\scriptstyle r\le j\le l-k-r-1\atop\scriptstyle i\le j}
\qst^{i-r-1}\ps_{j-i}(\qst)b^{s-j-k}\qst^{s-j-1} \cr
&\qquad\qquad+
\sum_{k+r+1\le i\le s} \;
\sum_{\scriptstyle l-k-r\le j\le l-k-1\atop\scriptstyle i\le j}
b^{i-r-k}\qst^{i-r-1}\ps_{j-i}(\qst)\qst^{s-j-1} \cr
&\qquad\qquad+
\sum_{k+r+1\le i\le s} \;
\sum_{\scriptstyle r\le j\le l-k-r-1\atop\scriptstyle i\le j}
b^{i-r-k}\qst^{i-r-1}\ps_{j-i}(\qst)b^{s-j-k}\qst^{s-j-1} \bigg). \cr
}$$
Here we have broken the sum into four parts, according to whether
$k+1\le i\le k+r$ or $k+r+1\le i\le s$, and independently according to whether 
$l-k-r\le j\le l-k-1$ or $r\le j\le l-k-r-1$.
(We note that the second and third double sums will vanish unless $c>5/2$,
and the fourth double sum will vanish unless $c>3$.)
The factor of $(H-1)^2$ accounts for the choice of $\rh(u)\in V\setminus\{0^k\}$ 
and $\si(u)\in W\setminus\{0^k\}$,
the factors preceding $\ps_{j-i}(\qst)$ in the summands account for the 
probability that the links of $u$ in ranks less than $i$ are idle,
the factors or $\ps_{j-i}(\qst)$ account for the probability that
the links of $u$ between $i$ and $j$ and not labelled $0^k$
are idle, and the factors following
$\ps_{j-i}(\qst)$ in the summands account for the probability that the links
of $u$ in ranks greater than $j$ are idle.
Bounding the factors $\ps_{j-i}(\qst)$ using Lemma 3.6, and bounding each
double summation by the number of terms (at most $(l-k)^2$) times
the largest term (which occurs for  $i=k+1$ and $j=l-k-r$ in the first
sum, and for $i=j$ for the remaining three sums), we obtain
$$\sigiv \le (H-1)^2 k (l-k)^2 \left(
\qst^k + 2b^{l-5k/2-1} \qst^{l-k-2} + b^{l-3k} \qst^{l-k-2} \right).$$
Thus we have
$$\eqalign{
{\sigiv\over\Ex[Y]} 
&\le k (l-k)^2 
\left({1\over (b\qst)^{l-2k}} + {2\over \qst b^{k/2+1}} + {1\over \qst b^k}
\right) \cr
&\to 0, \cr
}$$
since $b\qst>1$, $c>2$ and $b\ge 2$.
This completes the proof of the proposition, and with it the proof of Theorem
1.2.
\QED
\vfill\eject

\heading{4. Conclusion}

We have determined the limiting value of the linking probability
in spider-web networks with scale $k$ and depth $l$ when 
$l = ck$
with $c>1$.
The same method could be used when $l/k\to\infty$ but
$\(\log (l+1)\)/(k+1)\to 0$.
In this case, the phase transition would be less abrupt:
the limiting value of $Q$, and even its first derivative with respect to $q$,
would be continuous at the critical value $q_\infty = 1/b$, but the 
second derivative would be discontinuous.
There would be little to be gained by such networks, however,
over those with a large fixed value of $c$:
Their great cost would  decrease the
critical vacancy probability through only a small interval $[q_\infty, q_c]$,
and would provide only a small linking probability in this interval.

Another extension of our results would be to consider,
instead of the ``independent'' probability distribution on states
introduce by Lee [L1] and Le Gall [L2, L3], the ``coherent'' distribution
introduced by Pippenger [P1].
(The similar distribution introduced by Koverninski\u{\i} [K]
does not have an obvious generalization for $c>2$, and in any case it
does not seem likely that the additional independence in Koverninski\u{\i}'s
model would have much effect on its tractability for $c>2$.)

Yet another line of inquiry would be to consider the computational
complexity of path-search problems for spider-web networks with $c>2$,
using the link-probe model introduced by Lin and Pippenger [L4].
Some such results were obtained by Pippenger [P4] for $c=2$
(and these results are easily extended to the case $1<c<2$),
but even for $c=2$ the known results are far from definitive.
\sk

\heading{5. Acknowledgment}

The results reported in this paper were obtained 
during the fourth meeting of the Institute for Elementary Studies,
the Focused Research Group on Problems in Discrete Probability,
held 12--26 July 2003 at the Banff International Research Station
in Banff, Alberta, Canada.
The author is especially grateful to Yuval Peres, one of the organizers
of that meeting, for urging continued faith in the power of the
second-moment method.
\sk

\heading{6. References}

\ref B; S. R. Broadbent and J. M. Hammersley;
``Percolation Processes. I. Crystals and Mazes'';
Proceedings or the Cambridge Philosophical Society; 53 (1957) 629--641.

\ref C; F. R. K. Chung and F. K. Hwang;
``The Connection Pattern of Two Binary Trees'';
SIAM Journal on Algebraic and Discrete Methods; 1 (1980) 322--335.

\refbook G; G. Grimmett;
Percolation;
Springer-Verlag, New York, 1989
(Second edition: Springer-Verlag, Berlin, 1999).

\ref I; N. Ikeno;
``A Limit on Crosspoint Number'';
IEEE Transactions on Information Theory; 5 (1959) 187--196.

\ref K; I. V. Koverninski\u\i;
``Estimation of the Blocking Probability for Switching Circuits
by Means of Probability Graphs'';
Problems of Information Transmission; 11 (1975) 63--71.

\ref L1; C. Y. Lee;
``Analysis of Switching Networks'';
Bell System Technical Journal; 34 (1955) 1287--1315.

\ref L2; P. Le Gall; 
``\'Etude du blocage dans les syst\`emes de commutation 
t\'el\'ephoniques automatiques utilisant des commutateurs
\'electroniques du type crossbar'';
Annales des T\'el\'ecommunications; 11 (1956) 159--171, 180--194, 197.

\ref L3; P. Le Gall;
M\'ethod de calcul de l'encombrement dans les syst\`emes  
t\'el\'ephoniques automatiques a marquage;
Annales des T\'el\'ecommunications; 12 (1957) 374--386.

\ref L4; G. Lin and N. Pippenger;
``Routing Algorithms for Switching Networks with Probabilistic Traffic'';
Networks; 28 (1996) 21--29.

\ref P1; N. Pippenger;
``On Crossbar Switching Networks'';
IEEE Transactions on Communications; 23 (1975) 646--659.

\ref P2; N. Pippenger;
``The Blocking Probability of Spider-Web Networks'';
Random Structures \& Algorithms; 2 (1991) 121--149.

\ref P3; N. Pippenger;
``The Asymptotic Optimality of Spider-Web Networks'';
Discrete Applied Mathematics; 37/38 (1992) 437--450.

\ref P4; N. Pippenger;
``Upper and Lower Bounds for the Average-Case Complexity
of Path Search'';
Networks; 33 (1999) 249--259.

\ref T; K. Takagi;
``Design of Multi-Stage Link Systems by Means of
Optimal Channel Graphs'';
Electronic Communiations in Japan; 51A (1968) 37--46.

\bye